\documentclass[11pt,reqno]{amsart}
\usepackage{graphicx}
\usepackage{amscd}
\usepackage{amsmath}
\usepackage{amsfonts}
\usepackage{amssymb}
\newtheorem{theorem}{Theorem}[section]
\theoremstyle{plain}

\newtheorem{corollary}[theorem]{Corollary}
\newtheorem{proposition}[theorem]{Proposition}

\theoremstyle{definition}
\newtheorem{definition}[theorem]{Definition}
\newtheorem{example}[theorem]{Example}
\newtheorem{remark}[theorem]{Remark}
\newtheorem*{acknowledgements}{Acknowledgements}
\numberwithin{equation}{section}

\newcommand{\cala}{{\mathcal A}}

\newcommand{\cale}{{\mathcal E}}

\newcommand{\calp}{{\mathcal P}}

\newcommand{\frakh}{\mathfrak{h}}

\newcommand{\bleft}{[\![}
\newcommand{\bright}{]\!]}

\newcommand{\backl}{\mathbin{\vrule width1.5ex height.4pt\vrule height1.5ex}}

\begin{document}
\title[Leibniz Algebras, Courant Algebroids, and Multiplications]{Leibniz Algebras, Courant Algebroids, and Multiplications on Reductive
Homogeneous Spaces}
\author{Michael K. Kinyon}
\address{Department of Mathematics \& Computer Science\\
Indiana University\\
South Bend, IN\ 46634 USA}
\email{mkinyon@iusb.edu}
\urladdr{http://www.iusb.edu/\symbol{126}mkinyon}
\author{Alan Weinstein}
\address{Department of Mathematics\\
University of California\\
Berkeley, CA 94720 USA}
\email{alanw@math.berkeley.edu}
\urladdr{http://www.math.berkeley.edu/\symbol{126}alanw}
\thanks{Research partially supported by NSF\ grants DMS-96-25512 and DMS-99-71505 and
the Miller Institute for Basic Research.}
\thanks{This paper is in final form, and no version of it will be submitted for
publication elsewhere.}
\date{\today}

\begin{abstract}
We show that the skew-symmetrized product on every Leibniz algebra
$\mathcal{E}$ can be realized on a reductive complement to a subalgebra in a
Lie algebra. As a consequence, we construct a nonassociative multiplication on
$\mathcal{E}$ which, when $\mathcal{E}$ is a Lie algebra, is derived from the
integrated adjoint representation. We apply this construction to realize the
bracket operations on the sections of Courant algebroids and on the ``omni-Lie
algebras'' recently introduced by the second author.
\end{abstract}\maketitle

\section{Introduction}

\label{sec-intro}

Skew-symmetric bilinear operations which satisfy weakened versions of the
Jacobi identity arise from a number of constructions in algebra and
differential geometry. The purpose of this paper is to show how certain of
these operations, in particular the Courant brackets on the doubles of Lie
bialgebroids, can be realized in a natural way on the tangent spaces of
reductive homogeneous spaces. We use our construction to take steps toward
finding group-like objects which ``integrate'' these not-quite-Lie algebras.

The main ideas behind our construction come from work of K. Nomizu, K.
Yamaguti, and M. Kikkawa. Nomizu \cite{no:invariant} showed that affine
connections with parallel torsion and curvature are locally equivalent to
invariant connections on reductive homogeneous spaces, and that each such
space has a canonical connection for which parallel translation
along geodesics agrees with the natural action of the group. Yamaguti
\cite{ya:triple} characterized the torsion and curvature tensors of Nomizu's
canonical connection as pairs of algebraic operations, one bilinear and the
other trilinear, satisfying axioms defining what he called a ``general
Lie triple system'', and what Kikkawa later called a ``Lie triple
algebra''. In this paper, we will call these objects \emph{Lie-Yamaguti
algebras}. When the trilinear operation is zero, the bilinear operation is a
Lie algebra operation, and the homogeneous space is locally a Lie group on 
which the connection is the one which makes left-invariant vector fields 
parallel. Finally, Kikkawa \cite{ki:geometry} showed how to ``integrate'' 
Lie-Yamaguti algebras to nonassociative multiplications on reductive 
homogeneous spaces, and he characterized these multiplications axiomatically. 
Unfortunately, Kikkawa's construction when applied in our setting does
not quite reproduce the multiplication on a Lie group when the curvature 
is zero; rather it gives the loop operation 
$(x,y)\mapsto x+\exp(\mathrm{ad}\ x/2)y$ on the Lie algebra
itself. This limitation extends to our own work, so that we do not finally
succeed in finding the group-like object which we seek.

The starting point of our investigations was a skew-symmetric but non-Lie
bracket introduced by T.~Courant \cite{co:dirac}. It is defined on the direct
sum $\mathcal{E}=\mathcal{X}(P)\oplus\Omega^{1}(P)$ of the smooth vector
fields and $1$-forms on a differentiable manifold $P$ by
\begin{equation}
\lbrack\![(\xi_{1},\theta_{1}),(\xi_{2},\theta_{2})]\!\rbrack= \left(
[\xi_{1},\xi_{2}] , \mathcal{L}_{\xi_{1}}\theta_{2} - \mathcal{L}_{\xi_{2}%
}\theta_{1} - \frac{1}{2}d(i_{\xi_{1}}\theta_{2} - i_{\xi_{2}}\theta_{1})
\right)  , \label{eq:courant}%
\end{equation}
where $\mathcal{L}_{\xi}$ and $i_{\xi}$ are the operations of Lie derivative
and interior product by the vector field $\xi$. The term $-\frac{1}{2}%
d(i_{\xi_{1}}\theta_{2} - i_{\xi_{2}}\theta_{1}))$ will be especially
important to our discussion. It distinguishes the bracket from that on the
semidirect product of the vector fields acting on the 1-forms by Lie
derivation, and it spoils the Jacobi identity. On the other hand, with this
term (this was Courant's original motivation for introducing this bracket),
the graph of every Poisson structure $\Gamma\pi:\Omega^{1}(P)\to
\mathcal{X}(P)$ and every closed 2-form $\Gamma\omega:\mathcal{X}(P)\to
\Omega^{1}(P)$ is a subalgebra of $\mathcal{E}$. Thus, although $\mathcal{E}$
is not a Lie algebra, it contains many Lie algebras. Is there a ``group-like''
object associated to $\mathcal{E}$ which contains the (infinite-dimensional,
possibly local) Lie groups associated to the Lie subalgebras of $\mathcal{E}$?

When $P$ is the dual $V^{\ast}$ of a vector space $V$, truncating power series
expansions of vector fields and 1-forms leads to a bracket on the finite
dimensional space $\mathcal{E}=\mathfrak{gl}(V)\times V$. The bracket is
defined by
\begin{equation}
\lbrack\![(X,u),(Y,v)]\!]=\left(  [X,Y],\frac{1}{2}(Xv-Yu)\right)
\label{eq:omni}%
\end{equation}
for $X,Y\in\mathfrak{gl}(V)$ and $u,v\in V$. The factor of $1/2$ is now the
``spoiler'' of the Jacobi identity, but it yields the following nice property
of the bracket. If $\mathrm{ad}_{B}:V\rightarrow\mathfrak{gl}(V)$ is the
adjoint representation of any skew-symmetric operation $B:V\times V\rightarrow
V$, then $(V,B)$ is a Lie algebra if and only if the graph of $\mathrm{ad}%
_{B}$ is a subalgebra of $(\mathfrak{gl}(V)\times V,[\![\cdot,\cdot]\!])$,
which is then isomorphic to $(V,B)$ under the projection onto the second
factor \cite{we:omni}. 
The problem at the end of the last paragraph now becomes: is there a
group-like object associated to $(\mathfrak{gl}(V)\times V,[\![\cdot
,\cdot]\!])$ which contains Lie groups associated to all the Lie algebra
structures on $V$?

In this paper, we give a partial solution to these problems. We show that each
of the algebras denoted by $\mathcal{E}$ above can be embedded in a Lie
algebra $\mathcal{D}$ of roughly twice the size. There is a direct sum
decomposition $\mathcal{D} = \overline{\mathcal{E}}\oplus\mathcal{E}$
invariant under the adjoint action of $\overline{\mathcal{E}}$, a subalgebra.
The bracket on $\mathcal{E}$ is obtained from the bracket in $\mathcal{D}$ by
projection along $\overline{\mathcal{E}}$.

Denoting by $G$ and $H$ the groups associated to
$\mathcal{D}$ and $\overline{\mathcal{E}}$ respectively, we may identify
$\mathcal{E}$ with the tangent space at the basepoint of the reductive
homogeneous space $S(\mathcal{E})=G/H$. If $\mathcal{L}%
\subseteq\mathcal{E}$ is any Lie subalgebra, such as the graph of a Poisson
structure or of an adjoint representation, then there is a reductive
homogeneous space $S(\mathcal{L})$ which sits naturally inside $S(\mathcal{E}%
)$. This is as close as we have come so far to solving our problem. It is not
a complete solution since $S(\mathcal{L})$ is not a group. It does, however,
carry a multiplication which ``partially integrates'' the Lie algebra
structure on $\mathcal{L}$. The general procedure for constructing such
multiplications is as follows.

Any reductive homogeneous space $G/H$ corresponds to a reductive Lie algebra
decomposition $\mathfrak{g}=\mathfrak{h}\oplus\mathfrak{m}$; i.e.,
$\mathfrak{h}$ is a subalgebra with $[\mathfrak{h},\mathfrak{m}]\subseteq
\mathfrak{m}$. (Note that this is ``reductive'' in the sense of Nomizu
\cite{no:invariant}, not in the sense of having a completely reducible adjoint
representation.) In general $\mathfrak{m}$ is not a subalgebra, but a
neighborhood of the identity in $M=\exp(\mathfrak{m})\subseteq G$ can 
still be
identified with a neighborhood of the base point in $G/H$. Ignoring for
simplicity of exposition the restriction to neighborhoods of the identity
(which is in fact unnecessary in many cases), we may now take two elements $x$
and $y$ of $G/H$, multiply their representatives in $M$ to get a result in
$G$, and project to get a product $xy$ in $G/H$. Unless $M$ is a subgroup,
this multiplication will not, in general, be associative, but it will satisfy
the axioms which make $G/H$ into a group-like object called a
\emph{homogeneous Lie loop} \cite{ki:geometry}. Each Lie-Yamaguti subalgebra
$\mathfrak{l}$ of $\mathfrak{m}$ will then correspond to a subloop of $G/H$.

In the example $\mathcal{E}=\mathfrak{gl}(V)\times V$, the
Lie algebra $\mathcal{D}$ may be taken to be simply the semidirect product of
$\mathfrak{gl}(V)$ acting on the vector space $\mathfrak{gl}(V)\times V$ by
the direct sum of the adjoint and standard representations. The subalgebra
$\mathcal{E}$ is the first (i.e. the nonabelian) copy of $\mathfrak{gl}(V)$,
and the reductive complement is the graph of the mapping $(X,u)\mapsto X/2$
from the last two factors of $\mathfrak{gl}(V)\times\mathfrak{gl}(V)\times V$
to the first one. The construction for the original Courant bracket is
similar and is described in detail in \S4 below.

The construction of the enveloping Lie algebras $\mathcal{D}$ was in fact
worked out in a more general setting. It was observed in \cite{li-we-xu:manin}
that adding a symmetric term to the skew-symmetric Courant bracket
``improved'' some of its algebraic properties. Y.~Kosmann-Schwarzbach and
P.~Xu, as well as P.~\v{S}evera (all unpublished), observed that this
unskewsymmetrized operation, which we will denote by $\cdot$, satisfies the
derivation identity
\begin{equation*}
x\cdot(y\cdot z)=(x\cdot y)\cdot z+y\cdot(x\cdot z) 
\end{equation*}
defining what Loday \cite{lo:version} has called a \emph{Leibniz algebra}. The
central result in our paper (Theorems \ref{thm:envelope2} and
\ref{thm:recovery})
  is that the skew-symmetrization of every Leibniz
algebra structure can be extended in a natural way to a Lie-Yamaguti structure
and hence can be realized as the projection of a Lie algebra bracket onto a
reductive complement of a subalgebra. Thus, the Leibniz algebra is
``integrated'' to a homogeneous left loop. Whether this loop will
enable us to lasso Loday's elusive
``coquecigrue'' \cite{lo:version} remains to be seen.

%
%

In the last section of the paper, we discuss further directions for
research.  One is to attempt to take into account the fact that
the algebra $\mathcal{X} (P)\oplus\Omega^{1}(P)$ is also the space of
sections of a vector bundle, and to try to build a corresponding
structure into its enveloping Lie algebra.  A second goal is to
construct a group-like object for a Leibniz algebra which is actually
a group in the case of a Lie algebra. We end the paper with an idea
for constructing such an object as a quotient of a path space.

\begin{acknowledgements}
We would like to thank the many people from whom we have learned
important things about Courant algebroids, path spaces, and
nonassociative algebras, among them Anton Alekseev, Hans Duistermaat,
Nora Hopkins, Johannes Huebschmann, Atsushi Inoue, Yvette
Kos\-mann-Schwarzbach, Zhang-ju Liu, Kirill Mackenzie, Hala Pflugfelder,
Jon Phillips, Dmitry Roytenberg, Arthur Sagle, Pavol \v{S}evera, Jim
Stasheff, and Ping Xu. The second author would also like to thank
Setsuro Fujiie and Yoshitsugu Takei for their invitation to speak on
``Omni-Lie algebras'' and to publish a report in the proceedings of
the RIMS Workshop on Microlocal Analysis of Schr\"odinger
operators. Although that report had little to do with the subject of
the workshop, the stimulus to write a manuscript led to the posting of
a preprint \cite{we:omni} on the arXiv server and the subsequent
collaboration which has resulted in the present paper.
\end{acknowledgements}

\section{Leibniz algebras}

\label{sec-loday}

All vector spaces, algebras, etc. in this section will be over a ground field
$\mathbb{K}$ of characteristic $0$. Most results extend in obvious ways to
positive characteristic (not $2$), or even to commutative rings with unit. By
an \emph{algebra} $(\mathcal{E},\cdot)$, we will mean a vector space
$\mathcal{E}$ over $\mathbb{K}$ with a not necessarily associative bilinear
operation $\cdot:\mathcal{E}\times\mathcal{E}\rightarrow\mathcal{E}$. For
$x\in\mathcal{E}$, let $\lambda(x):\mathcal{E}\rightarrow\mathcal{E};y\mapsto
x\cdot y$ denote the left multiplication map. Let $\mathrm{Der}(\mathcal{E})$
denote the Lie subalgebra of $\mathfrak{gl}(\mathcal{E})$ consisting of the
derivations of $\mathcal{E}$. A linear map $\xi\in\mathfrak{gl}(\mathcal{E})$
is a derivation of $(\mathcal{E},\cdot)$ if and only if
\begin{equation}
\left[  \xi,\lambda(x)\right]  =\lambda(\xi x) \label{eq:lambda-equi}%
\end{equation}
for all $x\in\mathcal{E}$. For the class of algebras of interest to us, the
left multiplication maps have a stronger compatibility with the derivations.

\begin{definition}
\label{def:leibniz}An algebra $(\mathcal{E},\cdot)$ is called a \emph{Leibniz
algebra} if
\begin{equation}
x\cdot(y\cdot z)=(x\cdot y)\cdot z+y\cdot(x\cdot z) \label{eq:leibniz}%
\end{equation}
for all $x,y,z\in L$.
\end{definition}

Clearly, an algebra $(\mathcal{E},\cdot)$ is a Leibniz algebra if and only if
$\lambda(\mathcal{E})\subseteq\mathrm{Der}(\mathcal{E})$, or equivalently,
$\lambda:(\mathcal{E},\cdot)\to(\mathfrak{gl}(\mathcal{E}),[\cdot,\cdot])$ is
a homomorphism. Thus we have a homomorphism $\lambda:(\mathcal{E},\cdot
)\to(\mathrm{Der}(\mathcal{E}),[\cdot,\cdot])$ when $(\mathcal{E},\cdot)$ is a
Leibniz algebra.

Leibniz algebras were introduced by Loday \cite{lo:version}. (For this reason,
they have also been called ``Loday algebras'' \cite{ko:from}.) A
skew-symmetric Leibniz algebra structure is a Lie bracket; in this case,
(\ref{eq:leibniz}) is just the Jacobi identity. In particular, given a Leibniz
algebra $(\mathcal{E},\cdot)$, any subalgebra 
on which $\cdot$ is skew-symmetric is a Lie algebra, as is any
skew-symmetric quotient.

The skew-symmetrization of a Leibniz algebra $(\mathcal{E},\cdot)$ is an
interesting structure in its own right. We will denote the skew-symmetrized
operation by
\begin{equation}
\lbrack\![x,y]\!]=\frac{1}{2}\left(  x\cdot y-y\cdot x\right)  \label{eq:skew}%
\end{equation}
for $x,y\in\mathcal{E}$. In general, $(\mathcal{E},[\![\cdot,\cdot]\!])$ is
not a Lie algebra, i.e. $(\mathcal{E},\cdot)$ is not Lie-admissible
\cite{myung:malcev}. From (\ref{eq:leibniz}), we have that $\lambda
(x)\in\mathrm{Der}(\mathcal{E},[\![\cdot,\cdot]\!])$ for all $x\in\mathcal{E}%
$, and that $\lambda:(\mathcal{E},[\![\cdot,\cdot]\!])\to\mathrm{Der}%
(\mathcal{E})$ is a homomorphism of anticommutative algebras.

Let
\begin{equation*}
\mathcal{J}=\langle x\cdot x|x\in\mathcal{E}\rangle
\end{equation*}
be the two-sided ideal of $(\mathcal{E},\cdot)$ generated by all squares. Then
$\mathcal{J}$ contains all symmetric products $x\cdot y+y\cdot x$ for
$x,y\in\mathcal{E}$. Since $\lambda(x\cdot x)=\left[  \lambda(x),\lambda
(x)\right]  $ for all $x\in\mathcal{E}$, it follows that $\mathcal{J}%
\subseteq\ker(\lambda)$. Let $\mathcal{M}\subseteq\mathcal{E}$ be any ideal
containing $\mathcal{J}$. Since $x\cdot y+\mathcal{M}=-y\cdot x+(x\cdot
y+y\cdot x)+\mathcal{M}=-y\cdot x+\mathcal{M}$ for $x,y\in\mathcal{E}$, the
Leibniz product in $\mathcal{E}$ descends to a Lie bracket $[\cdot,\cdot]$ in
$\mathcal{E}/\mathcal{M}$. Conversely, if $\mathcal{M}\subseteq\mathcal{E}$ is
an ideal such that the quotient $\mathcal{E}/\mathcal{M}$ is a Lie algebra,
then for $x\in\mathcal{E}$, we have $x\cdot x+\mathcal{M}=\mathcal{M}$, and
thus $\mathcal{J}\subseteq\mathcal{M}$. In particular, $\mathcal{J}$ itself is
the smallest two-sided ideal of $\mathcal{E}$ such that
$\mathcal{E}/\mathcal{J}$ is a Lie algebra.

We now introduce one of our principal examples, which can be viewed as a
non-skew-symmetrized semidirect product of Lie algebras.

\begin{example}
\label{ex:hemi  &  demi}Let $(\mathfrak{h},[\cdot,\cdot])$ be a Lie algebra,
and let $V$ be an $\mathfrak{h}$-module with left action $\mathfrak{h}\times
V\rightarrow V:(\xi,x)\mapsto\xi x$. The induced left action of $\mathfrak{h}$
on $\mathfrak{h}\times V$ is just the restricted adjoint action of
$\mathfrak{h}$ on the semidirect product $\mathfrak{h}\ltimes V$:
\begin{equation}
\xi(\eta,y):=\left[  (\xi,0),(\eta,y)\right]  =([\xi,\eta],\xi y)
\label{eq:diag-action}%
\end{equation}
for $\xi,\eta\in\mathfrak{h}$, $y\in V$. Define a binary operation $\cdot$ on
$\mathcal{E}=\mathfrak{h}\times V$ by
\begin{equation*}
(\xi,x)\cdot(\eta,y)=\xi(\eta,y) 
\end{equation*}
for $\xi,\eta\in\mathfrak{h}$, $x,y\in V$; i.e.
\begin{equation}
(\xi,x)\cdot(\eta,y)=([\xi,\eta ],\xi y). \label{eq:hxV-product2}%
\end{equation}
Then $(\mathcal{E},\cdot)$ is a Leibniz algebra, and if $\mathfrak
{h}$ acts nontrivially on $V$, then $(\mathcal{E},\cdot)$ is not a Lie
algebra. We call $\mathcal{E}$ with this Leibniz algebra structure the
\emph{hemisemidirect product} of $\mathfrak{h}$ with $V$, and denote it by
$\mathfrak{h}\ltimes_{H}V$.

The skew-symmetrized product in $(\mathcal{E},\cdot)$ is
\begin{equation}
\lbrack\![(\xi,x),(\eta,y)]\!]=\left(  [\xi,\eta],\frac{1}{2}(\xi y-\eta
x)\right)  \label{eq:demi}%
\end{equation}
for $\xi,\eta\in\mathfrak{h}$, $x,y\in V$. We call $\mathcal{E}$ with the
bracket $[\![\cdot,\cdot]\!]$ the \emph{demisemidirect product} of
$\mathfrak{h}$ with $V$, and we enote it by $\mathfrak
{h}\ltimes_{D}V$ As we noted in \S\ref{sec-intro}, the factor of 1/2 generally
spoils the Jacobi identity for this bracket. We have
\begin{align*}
\mathcal{J}  &  =\left\{  0\right\}  \times\mathfrak{h}V\\
\ker(\lambda)  &  =\left\{  \xi\in\mathfrak{z}(\mathfrak{h})|\xi V=0\right\}
\times V
\end{align*}
where $\mathfrak{z}(\mathfrak{h})$ is the center of $\mathfrak{h}$. If the
representation of $\mathfrak{h}$ on $V$ is faithful and if $\mathfrak{h}V=V$,
then $\mathcal{J}=\ker(\lambda)$. For example, if $\mathfrak{h}=\mathfrak
{gl}(V)$, then $\mathcal{J}=\ker(\lambda)=\left\{  0\right\}  \times V$.

Let $\pi_{\mathfrak{h}}:\mathcal{E}\rightarrow\mathfrak{h}$ denote the
projection onto the first factor. Then $\pi_{\mathfrak{h}}$ is $\mathfrak{h}%
$-equivariant, i.e.
\begin{equation}
\left[  \xi,\pi_{\mathfrak{h}}(\eta,x)\right]  =\pi_{\mathfrak{h}}\left(
\xi(\eta,x)\right)  \label{eq:pi-equi}%
\end{equation}
for all $\xi,\eta\in\mathfrak{h}$, $x\in V$. The homomorphism $\lambda
:\mathcal{E}\rightarrow\mathrm{Der}(\mathcal{E})$ factors through
$\pi_{\mathfrak{h}}$ and the action (\ref{eq:diag-action}):
\begin{equation}
\lambda(\xi,x)(\eta,y)=\pi_{\mathfrak{h}}(\xi,x)(\eta,y). \label{eq:pi-action}%
\end{equation}
\end{example}

The preceding paragraph motivates the following definition.

\begin{definition}
\label{def:envelope}Let $(\mathcal{E},\cdot)$ be a Leibniz algebra. Let
$\mathfrak{h}$ be a Lie algebra with a derivation action $\mathfrak
{h}\to\mathrm{Der}(\mathcal{E})$, i.e. $\xi(x\cdot y)=\xi x\cdot y+x\cdot\xi
y$ for all $\xi\in\mathfrak{h}$, $x,y\in\mathcal{E}$. Let $f:\mathcal{E}%
\to\mathfrak{h}$ be an $\mathfrak{h}$-equivariant linear map such that the
diagram
\begin{equation}%
\begin{array}
[c]{ccc}%
\mathcal{E} & \overset{\lambda}{\longrightarrow} & \mathrm{Der}(\mathcal{E})\\
f\downarrow & \nearrow & \\
\mathfrak{h} &  &
\end{array}
\label{eq:diagram}%
\end{equation}
commutes. Let $\mathfrak{g}=\mathfrak{h}\ltimes\mathcal{E}$ be the semidirect
product Lie algebra, $\mathcal{E}$ being considered as a Lie algebra with the
zero bracket. We say that the triple $(\mathfrak{g},\mathfrak{h},f)$ is an
\emph{enveloping Lie algebra} of $(\mathcal{E},\cdot)$. We will justify this
name with Theorem \ref{thm:recovery} below.
\end{definition}

In terms of equations, the $\mathfrak{h}$-equivariance of $f$ and the
commuting of (\ref{eq:diagram}) are expressed by
\begin{align}
\left[  \xi,f(x)\right]   &  =f(\xi x)\label{eq:h-equi}\\
f(x)y  &  =\lambda(x)y \label{eq:f(x)y=x.y}%
\end{align}
for all $x,y\in\mathcal{E}$, $\xi\in\mathfrak{h}$, $x\in\mathcal{E}$. (The
$\mathfrak{h}$-equivariance of $f$ implies that
$f(\mathcal{E})$ is a Lie ideal of $\mathfrak{h}$.) Properties (\ref{eq:h-equi}%
) and (\ref{eq:f(x)y=x.y}) imply that
\begin{equation*}
f(x\cdot y)=f(f(x)y)=\left[  f(x),f(y)\right]  
\end{equation*}
for $x,y\in\mathcal{E}$; that is, $f$ is a homomorphism of Leibniz algebras.
In particular, $f(x\cdot x)=0$ for all $x\in\mathcal{E}$, and
\begin{equation*}
f \left(  [\![x,y]\!]\right)  = [f(x),f(y)]. 
\end{equation*}
By (\ref{eq:f(x)y=x.y}), if $f(x)=0$, then $\lambda(x)=0$. Therefore we have
the inclusions
\begin{equation*}
\mathcal{J}\subseteq\ker(f)\subseteq\ker(\lambda). 
\end{equation*}

Just as (\ref{eq:pi-equi}) and (\ref{eq:pi-action}) motivated the definition
of enveloping Lie algebra, so do they imply the following.

\begin{proposition}
\label{prop:hemi-env}$(\mathfrak{h}\ltimes\mathcal{E},\mathfrak{h}%
,\pi_{\mathfrak{h}})$ is an enveloping Lie algebra for the hemi\-semi\-direct
product $\mathcal{E}=\mathfrak{h}\ltimes_{H}V$.
\end{proposition}

Every Leibniz algebra $(\mathcal{E},\cdot)$ has enveloping Lie algebras. By
(\ref{eq:leibniz}), $\lambda(\mathcal{E})$ is a Lie subalgebra of
$\mathrm{Der}(\mathcal{E})$. 
Since left multiplication maps are
$\mathrm{Der}(\mathcal{E})$-equivariant (see (\ref{eq:lambda-equi})), the
following holds.

\begin{proposition}
\label{prop:envelope1}Let $(\mathcal{E},\cdot)$ be a Leibniz algebra,
and let $\mathfrak{h}$ be a Lie algebra satisfying
$\lambda(\mathcal{E}) \subseteq \mathfrak{h}\subseteq
\mathrm{Der}(\mathcal{E})$. Then $(\mathfrak
{h}\ltimes\mathcal{E},\mathfrak{h},\lambda)$ is an enveloping Lie algebra for
$(\mathcal{E},\cdot)$.
\end{proposition}

If $(\mathfrak{g},\mathfrak{h},f)$ is an enveloping Lie algebra of a Leibniz
algebra $(\mathcal{E},\cdot)$, then so is $(\tilde{\mathfrak{g}}%
,f(\mathcal{E}),f)$ where $\tilde{\mathfrak{g}}=f(\mathcal{E})\ltimes
\mathcal{E}$. For our purposes, the case where $f(\mathcal{E})=\mathfrak{h}$
is of most interest. In this case we have
\[
\mathfrak{h}\cong\mathcal{E}/\ker(f).
\]
Conversely, let $\mathcal{M}$ be an ideal of $(\mathcal{E},\cdot)$ satisfying
$\mathcal{J}\subseteq\mathcal{M}\subseteq\ker(\lambda)$. Then $\mathfrak
{h}:=\mathcal{E}/\mathcal{M}$ is a Lie algebra, and the Leibniz algebra
homomorphism $\lambda:\mathcal{E}\rightarrow\mathrm{Der}(\mathcal{E})$
descends to a Lie algebra homomorphism $\overline{\lambda}:\mathfrak
{h}\rightarrow\mathrm{Der}(\mathcal{E})$ such that
\[%
\begin{array}
[c]{ccc}%
\mathcal{E} & \overset{\lambda}{\longrightarrow} & \mathrm{Der}(\mathcal{E})\\
q\downarrow & \underset{\bar{\lambda}}{\nearrow} & \\
\mathfrak{h} &  &
\end{array}
\]
commutes, where $q:\mathcal{E}\rightarrow\mathfrak{h}$ is the natural
projection. The following result, which is our main construction of enveloping
Lie algebras, is an immediate consequence of these considerations.

\begin{theorem}
\label{thm:envelope2}Let $(\mathcal{E},\cdot)$ be a Leibniz algebra, let
$\mathcal{M}$ be an ideal of $\mathcal{E}$ satisfying $\mathcal{J}%
\subseteq\mathcal{M}\subseteq\ker(\lambda)$, and let $\mathfrak{h}%
=\mathcal{E}/\mathcal{M}$. Then $(\mathfrak{h}\ltimes\mathcal{E},\mathfrak
{h},q)$ is an enveloping Lie algebra for $(\mathcal{E},\cdot)$.
\end{theorem}

\begin{example}
\label{ex:hemi-env}Let $\mathcal{E}=\mathfrak{h}\ltimes_{H}V$ be the
hemisemidirect product, and assume that the representation of $\mathfrak{h}$
on $V$ is faithful and that $\mathfrak{h}V=V$. Then $\mathcal{J}=\ker
(\lambda)$, and $\mathfrak{h}\cong\mathcal{E}/\mathcal{J}$. We may identify
the natural projection $q:\mathcal{E}\rightarrow\mathcal{E}/\mathcal{J}$ with
the projection $\pi_{\mathfrak
{h}}:\mathcal{E}\rightarrow\mathfrak{h}$ onto the first factor. In this case,
the enveloping Lie algebra of Proposition \ref{prop:hemi-env} is that of
Theorem \ref{thm:envelope2}.
\end{example}

\begin{remark}
Let $\mathfrak{h}$ be a Lie algebra, let $\mathcal{E}$ be a left $\mathfrak
{h}$-module, and let $f:\mathcal{E}\rightarrow\mathfrak{h}$ be a $\mathfrak
{h}$-equivariant linear map. Define a binary operation $\cdot$ on
$\mathcal{E}$ by $x\cdot y=f(x)y$. Then $(\mathcal{E},\cdot)$ is clearly a
Leibniz algebra. We will discuss this point further in \S3. Loday and
Pirashvili \cite{lo:tensor} have shown that $f:\mathcal{E}\rightarrow
\mathfrak{h}$ can be considered to be a Lie algebra object in what they call
the infinitesimal tensor category $\mathcal{LM}$ of linear mappings. The
assignments $(f:\mathcal{E}\rightarrow\mathfrak
{h})\rightsquigarrow(\mathcal{E},\cdot)$ and $(\mathcal{E},\cdot
)\rightsquigarrow(q:\mathcal{E}\rightarrow\mathcal{E}/\mathcal{J})$ are
adjoint functors between the category of Lie algebra objects in $\mathcal{LM}$
and the category of Leibniz algebras.
\end{remark}

We now move on to our main result, which is
to show how to recover the skew-symmetrized structure $(\mathcal{E}%
,[\![\cdot,\cdot]\!])$
of a Leibniz algebra $(\mathcal{E},\cdot)$ from the Lie algebra
structure of an
enveloping Lie algebra
$(\mathfrak{g},\mathfrak{h},f)$.  Let
$\pi_{\mathcal{E}}:\mathfrak{g}\rightarrow\mathcal{E}$ denote the projection
onto the second factor. For each $s\in\mathbb{K}$, define a section
$\sigma_{s}:\mathcal{E}\rightarrow\mathfrak{g}$\ of $\pi_{\mathcal{E}}$ by
\begin{equation}
\sigma_{s}(x)=(sf(x),x) \label{eq:section}%
\end{equation}
for $x\in\mathcal{E}$. The image 
\begin{equation*}
\cale_s = \sigma_{s}\left(  \mathcal{E}\right)  =\left\{  (sf(x),x)|x\in\mathcal{E}%
\right\}  
\end{equation*}
is a copy of $\mathcal{E}$ which is a complement of $\mathfrak{h}=\ker
(\pi_{\mathcal{E}})$. We will write the corresponding vector space
decomposition of $\mathfrak{g}$ as
\begin{equation}
\mathfrak{g}\cong\mathfrak{h}\oplus\mathcal{E}_s. \label{eq:g-decomp}%
\end{equation}
Note that the case $s=0$ is just the semidirect product of $\mathfrak{h}$ with
$\mathcal{E}$.

Since $\mathcal{E}_s$ is essentially the graph of the $\mathfrak{h}%
$-equivariant map $sf:\mathcal{E}\rightarrow\mathfrak{h}$, $\sigma_{s}$ itself
is equivariant for the adjoint action (\ref{eq:diag-action}) of $\mathfrak{h}$
on $\mathfrak{g}$. Indeed, for $x\in\mathcal{E}$, $\xi\in\mathfrak{h}$, we
have
\begin{align*}
\xi\ \sigma_{s}(x)  &  =\left(  [\xi,sf(x)],\xi x\right) \\
&  =(sf(\xi x),\xi x)\\
&  =\sigma_{s}(\xi x).
\end{align*}
This shows that the bracket of an element of $\mathfrak{h}$ with an element of
$\mathcal{E}$ relative to the decomposition (\ref{eq:g-decomp}) agrees with
the action of $\mathfrak{h}$ on $\mathcal{E}$, independently of the value of
$s$:
\[
\left[  \xi,x\right]  =\xi x.
\]
Since
\begin{equation}
\left[  \mathfrak{h},\mathcal{E}\right]  \subseteq\mathcal{E}
\label{eq:E-reductive}%
\end{equation}
the decomposition (\ref{eq:g-decomp}) is \emph{reductive} \cite{no:invariant}.

Now we consider the bracket of two elements of
$\mathcal{E}\cong\mathcal{E}_s$  relative to (\ref{eq:g-decomp}). Here, unlike
(\ref{eq:E-reductive}), we expect the result to depend on $s$. For
$x,y\in\mathcal{E}$, we have
\begin{align*}
\left[  \sigma_{s}(x),\sigma_{s}(y)\right]   &  =\left[
(sf(x),x),(sf(y),y)\right] \\
&  =(s^{2}[f(x),f(y)],s(f(x)y-f(y)x))\\
&  =(s^{2}f([\![x,y]\!]),2s[\![x,y]\!])\\
&  =(-s^{2}f([\![x,y]\!]),0)+(2s^{2}f([\![x,y]\!]),2s[\![x,y]\!])\\
&  =(-s^{2}f([\![x,y]\!]),0)+\sigma_{s}(2s[\![x,y]\!]).
\end{align*}
We use the components of the result of this computation to define
skew-symmetric bilinear maps $[\![\cdot,\cdot]\!]_{s}:\mathcal{E}%
\times\mathcal{E}\rightarrow\mathcal{E}$ and $\Delta_{s}:\mathcal{E}%
\times\mathcal{E}\rightarrow\mathfrak{h}$ by
\begin{align}
\lbrack\![x,y]\!]_{s}  &  =\pi_{\mathcal{E}}\left(  [\sigma_{s}(x),\sigma
_{s}(y)]\right)  =2s[\![x,y]\!]\label{eq:proj-bracket}\\
\Delta_{s}(x,y)  &  =\pi_{\mathfrak{h}}\left(  [\sigma_{s}(x),\sigma
_{s}(y)]\right)  =-s^{2}f([\![x,y]\!]) \label{eq:proj-delta}%
\end{align}
for $x,y\in\mathcal{E}$.

Observe that the choice $s=1/2$ recovers the original skew-symmetrized Leibniz
product on $\mathcal{E}$:
\begin{equation}
\lbrack\![x,y]\!]_{1/2}=[\![x,y]\!] \label{eq:proj=original}%
\end{equation}
for $x,y\in\mathcal{E}$. This gives us the following result.

\begin{theorem}
\label{thm:recovery}Let $(\mathcal{E},\cdot)$ be a Leibniz algebra with
enveloping Lie algebra $(\mathfrak{g},\mathfrak{h},f)$, and let $\sigma
_{1/2}:\mathcal{E}\rightarrow\mathfrak{g}$ be defined by (\ref{eq:section}).
Then the skew-symmetrized product (\ref{eq:skew}) is given in terms of the Lie
bracket in $\mathfrak{g}$ by
\[
\lbrack\![x,y]\!]=\pi_{\mathcal{E}}\left(  \left[  \sigma_{1/2}(x),\sigma
_{1/2}(y)\right]  \right)
\]
for $x,y\in\mathcal{E}$.
\end{theorem}

In \S5, we will show that the map $\Delta_{s}:\mathcal{E}\times\mathcal{E}%
\rightarrow\mathfrak{h}$ defined by (\ref{eq:proj-delta}) introduces
additional structure into $\mathcal{E}$.

\begin{example}
\label{ex:hemi-recovery}Let $\mathcal{E}=\mathfrak{h}\ltimes_{H}V$ be the
hemisemidirect product, and let $(\mathfrak{g},\mathfrak{h},\pi_{\mathfrak{h}%
})$ be the enveloping Lie algebra obtained in Proposition \ref{prop:hemi-env}.
The section $\sigma_{s}:\mathcal{E}\rightarrow\mathfrak{g}$ is given by
\[
\sigma_{s}(\xi,x)=\left(  s\xi,\xi,x\right)
\]
for $\xi\in\mathfrak{h}$, $x\in V$. A calculation shows that the skew
symmetric maps $[\![\cdot,\cdot]\!]_{s}:\mathcal{E}\times\mathcal{E}%
\rightarrow\mathcal{E}$ and $\Delta_{s}:\mathcal{E}\times\mathcal{E}%
\rightarrow\mathfrak{h}$ are given by
\begin{align}
\lbrack\![(\xi,x),(\eta,y)]\!]_{s}  &  =2s\left(  \left[  \xi,\eta\right]
,\frac{1}{2}(\xi y-\eta x)\right) \label{eq:demi-proj}\\
\Delta_{s}\left(  (\xi,x),(\eta,y)\right)   &  =-s^{2}\left[  \xi,\eta\right]
\label{eq:demi-delta}%
\end{align}
for $\xi,\eta\in\mathfrak{h}$, $x,y\in V$. If we take $s=1/2$, (\ref{eq:demi-proj})
reduces to the demisemidirect product bracket (\ref{eq:demi}).
\end{example}

\section{Omni-Lie and Omni-Leibniz Algebras}

\label{sec-omni}

In this section we will show that \emph{every} Leibniz algebra can be embedded
in a hemisemidirect product Leibniz algebra. 

Let $(\mathcal{E},\cdot)$ be a Leibniz algebra with enveloping Lie algebra
$(\mathfrak{g},\mathfrak{h},f)$, and let $\sigma_{s}:\mathcal{E}%
\rightarrow\mathfrak{g}$ be the section of $\pi_{\mathcal{E}}$ defined by
(\ref{eq:section}). The vector space $\mathfrak{g}=\mathfrak{h}\times
\mathcal{E}$ has (at least) four natural algebra structures which are related
to the structure of $\mathcal{E}$. First, $\mathfrak{g}=\mathfrak{h}%
\ltimes\mathcal{E}$ is a semidirect product Lie algebra. We have shown in
Theorem \ref{thm:recovery} that the section $\sigma_{1/2}$ can be used to
recover the skew-symmetrized structure $(\mathcal{E},[\![\cdot,\cdot]\!])$
from the semidirect product $\mathfrak{g}=\mathfrak{h}\ltimes\mathcal{E}$.
Second, $\mathfrak{g}=\mathfrak{h}\times\mathcal{E}$ also has the structure of
a direct product of Leibniz algebras. As noted before, each 
$\mathcal{E}_s$ is the graph of the map $sf:\mathcal{E}\rightarrow
\mathfrak{h}$. Since each $sf$ is a homomorphism of Leibniz algebras, each
$\mathcal{E}_s$ is trivially a subalgebra of the direct product
$\mathfrak{g}=\mathfrak{h}\times\mathcal{E}$, which, for $s\neq0$, is
isomorphic to $(\mathcal{E},\cdot)$ under the map $(1/s)\pi_{\mathcal{E}%
}|_{\mathcal{E}_s}$.

Third and fourth, $\mathfrak{g}$ also has the hemisemidirect product structure
$\mathfrak{g}=\mathfrak{h}\ltimes_{H}\mathcal{E}$ and the demisemidirect
product structure $\mathfrak{g}=\mathfrak{h}\ltimes_{D}\mathcal{E}$. The next
result shows how these are related to the Leibniz and skew-symmetrized Leibniz
algebra structures, respectively, on $\mathcal{E}$.

\begin{proposition}
\label{prop:subalgebras}
\begin{enumerate}
\item $\sigma_{1}:(\mathcal{E},\cdot)\rightarrow\mathfrak{h}\ltimes
_{H}\mathcal{E}$ is a monomorphism of Leibniz algebras.

\item $\sigma_{1}:(\mathcal{E},[\![\cdot,\cdot]\!])\rightarrow\mathfrak
{h}\ltimes_{D}\mathcal{E}$ is a monomorphism of skew-symmetrized Leibniz algebras.
\end{enumerate}
\end{proposition}

\begin{proof}
For $x,y\in\mathcal{E}$, we compute
\begin{align*}
\sigma_{1}(x)\cdot\sigma_{1}(y)  &  =\left(  f(x),x\right)  \cdot\left(
f(y),y\right) \\
&  =\left(  \left[  f(x),f(y)\right]  ,f(x)y\right) \\
&  =\left(  f(x\cdot y),x\cdot y\right) \\
&  =\sigma_{1}(x\cdot y)
\end{align*}
where in the penultimate equality, we have used the fact that $f$ is a
homomorphism of Leibniz algebras and (\ref{eq:f(x)y=x.y}). This establishes
(1), and (2) follows from (1).
\end{proof}

Since every Leibniz algebra has enveloping Lie algebras, we have the following.

\begin{corollary}
\begin{enumerate}
\item  Every Leibniz algebra can be embedded as a subalgebra in a
hemisemidirect product.

\item  Every skew-symmetrized Leibniz algebra can be embedded as a subalgebra
in a demisemidirect product.
\end{enumerate}
\end{corollary}

By Propositions \ref{prop:envelope1} and \ref{prop:subalgebras}, the
monomorphism $x\mapsto(\lambda(x),x)$ embeds a given Leibniz algebra
$(\mathcal{E},\cdot)$ as a subalgebra of the hemisemidirect product
$\mathrm{Der}(\mathcal{E},\cdot)\ltimes_{H}\mathcal{E}$, which in turn can be
embedded as a subalgebra of the hemisemidirect product $\mathfrak
{gl}(\mathcal{E})\ltimes_{H}\mathcal{E}$. While different Leibniz algebra
structures on the vector space $\mathcal{E}$ can lead to different derivation
algebras, $\mathfrak{gl}(\mathcal{E})\ltimes_{H}\mathcal{E}$ contains
\emph{all} Leibniz algebra structures on $\mathcal{E}$. We now show that, in
fact, this exactly characterizes the Leibniz algebras among all algebra
structures on $\mathcal{E}$.

Let $(\mathcal{E},\cdot)$ be an algebra, and let%
\begin{equation*}
\mathcal{G}_{\lambda}=\left\{  (\lambda(x),x)|x\in\mathcal{E}\right\}
\end{equation*}
denote the graph of $\lambda$ as a subspace of $\mathfrak{gl}(\mathcal{E}%
)\times\mathcal{E}$.

\begin{proposition}
\label{prop:leibniz-graph}An algebra $(\mathcal{E},\cdot)$ is a Leibniz
algebra if and only if $\mathcal{G}_{\lambda}$ is a subalgebra of
$\mathfrak{gl}(\mathcal{E})\ltimes_{H}\mathcal{E}$. If these conditions hold,
then $\pi_{\mathcal{E}}|_{\mathcal{G}_{\lambda}}$ is an isomorphism from
$(\mathcal{G}_{\lambda},\cdot)$ to $(\mathcal{E},\cdot)$.
\end{proposition}

\begin{proof}
For $x,y\in V$, we have
\begin{equation}
(\lambda(x),x)\cdot(\lambda(y),y)=([\lambda(x),\lambda(y)],x\cdot y).
\label{eq:graph-product}%
\end{equation}
Thus $\mathcal{G}_{\lambda}$ is closed under the product $\cdot$ if and only
if $\lambda$ is a homomorphism. The remaining assertions are clear.
\end{proof}

We see from Proposition \ref{prop:leibniz-graph} that the class of all Leibniz
algebra structures on $\mathcal{E}$ corresponds to the class of all linear
maps from $\mathcal{E}$ to $\mathfrak{gl}(\mathcal{E})$ whose graphs are
subalgebras of the hemisemidirect product. Recalling that a skew-symmetric
subalgebra of a Leibniz algebra is a Lie algebra, we see that $\mathfrak
{gl}(\mathcal{E})\ltimes_{H}\mathcal{E}$ also contains all Lie algebra
structures on $\mathcal{E}$. From (\ref{eq:graph-product}), it is immediate
that the operation $\cdot$ in $\mathcal{G}_{\lambda}$ is skew-symmetric if and
only if the operation $\cdot$ in $\mathcal{E}$ is skew-symmetric.

\begin{corollary}
\label{coro:lie-equiv}An algebra $(\mathcal{E},\cdot)$ is a Lie algebra if and
only if $\mathcal{G}_{\lambda}$ is a Lie subalgebra of $\mathfrak
{gl}(\mathcal{E})\ltimes_{H}\mathcal{E}$. If these conditions hold, then
$\pi_{\mathcal{E}}|_{\mathcal{G}_{\lambda}}$ is an isomorphism from
$(\mathcal{G}_{\lambda},\cdot)$ to $(\mathcal{E},\cdot)$.
\end{corollary}

In case $\mathcal{E}=\mathbb{R}^{n}$, the demisemidirect product
$(\mathfrak{gl}(\mathcal{E})\ltimes_{D}\mathcal{E},[\![\cdot,\cdot]\!])$ is
the ``omni-Lie algebra'' with bracket (\ref{eq:omni})
of \cite{we:omni}. Symmetrizing the Leibniz product
in $\mathfrak{gl}(\mathcal{E})\ltimes_{H}\mathcal{E}$ defines a commutative
product $\circ$ by
\[
(\xi,x)\circ(\eta,y)=\left(  0,\frac{1}{2}\left(  \xi y+\eta x\right)
\right)  ,
\]
which was called an ``$\mathcal{E}$-valued bilinear form'' in \cite{we:omni}.
The hemisemidirect product Leibniz algebra
$\mathfrak{gl}(\mathcal{E})\ltimes_{H}\mathcal{E}$
combines both of these structures. The following restatement of Corollary
\ref{coro:lie-equiv} encompasses both Proposition 1 in \cite{we:omni} and the
subsequent discussion.

\begin{corollary}
\label{coro:lie-equiv2}An algebra $(\mathcal{E},\cdot)$ is a Lie algebra if
and only if (i) for all $(\xi,x),(\eta,y)\in\mathcal{G}_{\lambda}$,
$(\xi,x)\circ(\eta,y)=0$, and (ii) $(\mathcal{G}_{\lambda},[\![\cdot
,\cdot]\!])$ is a Lie subalgebra of $\mathfrak
{gl}(\mathcal{E})\ltimes_{D}\mathcal{E}$.
\end{corollary}

In particular, if $(\mathcal{E},\cdot)$ is an anticommutative algebra, then it
is a Lie algebra if and only if the graph of its adjoint representation is a
subalgebra of the demisemidirect product $\mathfrak
{gl}(\mathcal{E})\ltimes_{D}\mathcal{E}$.

\section{Courant Algebroids}

\label{sec-courant}

In this section, we will apply the methods of \S\ref{sec-loday} to Courant
algebroids, which include as special cases the doubles of Lie bialgebras and
the bundles $TP\oplus T^{\ast}P$ with the bracket on sections given by
(\ref{eq:courant}). The following definition was introduced in
\cite{li-we-xu:manin}.

\begin{definition}
\label{def:courant} A \emph{Courant algebroid} is a vector bundle $E\to P$
equipped with a nondegenerate symmetric bilinear form
$\langle\!\cdot,\cdot \!\rangle$ on
the bundle, a skew-symmetric bracket $[\![\cdot,\cdot]\!]$ on the space
$\mathcal{E}=\Gamma(E)$ of smooth sections of $E$, and a bundle map $\rho:E\to
TP$ such that, for any $x,y,z\in\mathcal{E}$ and $f,g\in\mathcal{A}=C^{\infty
}(P)$:

\begin{enumerate}
\item $\sum_{(x,y,z)}[\![ [\![ x,y]\!],z]\!]=\mathcal{D} T(x,y,z);$

\item $\rho[\![ x,y]\!]=[\rho x, \rho y];$

\item $[\![ x,fy]\!]=f[\![ x,y]\!]+(\rho(x)f)y-\langle x, y\rangle\mathcal{D}
f ;$

\item $\rho\mbox{\tiny{$\circ$}} \mathcal{D} =0$; i.e., for any $f, g $,
$\langle\mathcal{D} f, \mathcal{D} g\rangle=0$;

\item $\rho(x) \langle y,z\rangle=\langle[\![ x , y]\!]+\mathcal{D} \langle x
,y\rangle, z \rangle+\langle y, [\![ x , z]\!]+\mathcal{D} \langle x ,z
\rangle\rangle$,
\end{enumerate}
Here, $\sum_{(x,y,z)}$ denotes the sum over cyclic permutations of $x$,
$y$, and $z$, $T(x,y,z)$ is the function on $P$ defined by:
\begin{equation*}
T(x,y,z)=\frac{1}{3}\sum_{(x,y,z)}\langle\lbrack\![x,y]\!],z\rangle
\end{equation*}
and $\mathcal{D}:\mathcal{A}\rightarrow\mathcal{E}$ is the map defined by
$\mathcal{D}=\frac{1}{2}\beta^{-1}\rho^{\ast}d$, where $\beta$ is the
isomorphism between $E$ and $E^{\ast}$ given by the bilinear form. In other
words,
\begin{equation*}
\langle\mathcal{D}f,x\rangle=\frac{1}{2}\rho(x)f. 
\end{equation*}
\end{definition}

It was already noted in \cite{li-we-xu:manin} that adding the symmetric term
$(x,y)\mapsto\mathcal{D}\langle x,y\rangle$ to the bracket of a Courant
algebroid leads to an operation $\cdot$ with nicer properties. In fact, as
noted independently by Kosmann-Schwarzbach, \v{S}evera, and Xu, the operation
$\cdot$ makes $\mathcal{E}$ into a Leibniz algebra. (Details may be found in
\S2.6 of \cite{ro:courant}.)

This time, we will use our main construction of Theorem \ref{thm:envelope2} to
find an enveloping Lie algebra for $\mathcal{E}$. By the definition of the
operation $\cdot$, the image $\mathcal{D}(\mathcal{A})$ contains the ideal
$\mathcal{J}$ generated by squares. In fact, by the identity $x\cdot
\mathcal{D}f=\mathcal{D}\langle x,\mathcal{D}f\rangle$ (Lemma 2.6.2 of
\cite{ro:courant}), $\mathcal{D}(\mathcal{A})$ is an ideal in $\mathcal{E}$,
so it can play the role of $\mathcal{M}$ in our general construction; i.e.
$\mathcal{E}/\mathcal{D}(\mathcal{A})$ is a Lie algebra acting on
$\mathcal{E}$. We may therefore form the semidirect product Lie algebra
$\mathfrak{g}=\mathcal{E}/\mathcal{D}(\mathcal{A})\ltimes\mathcal{E}$, which
functions as an enveloping Lie algebra. Continuing with the application of our
general construction in Theorem \ref{thm:recovery}, we see that the Courant
algebroid bracket on $\mathcal{E}$ is obtained from the Lie algebra bracket on
$\mathfrak{g}=\mathcal{E}/\mathcal{D}(\mathcal{A})\times\mathcal{E}$ by
identification of $\mathcal{E}$ with the graph of $\frac{1}{2}$ times the
quotient map from $\mathcal{E}$ to $\mathcal{E}/\mathcal{D}(\mathcal{A})$, by
projection along the subalgebra $\mathcal{E}/\mathcal{D}(\mathcal{A})\times\{0\}$.

\begin{example}
\label{ex:point}
If $P$ is a single point, $\mathcal{E}$ is just a Lie algebra with an
invariant symmetric bilinear form, and $\mathcal{D}=0$.  The
enveloping Lie algebra is then $\mathcal{E} \ltimes \mathcal{E}$.  The
bracket on $\mathcal{E}$ is recovered by projection along the first
factor onto the subspace (not a subalgebra!) $\{(\frac{1}{2}x,x)|x\in
\mathcal{E}\}.$
\end{example}

For the original Courant bracket (\ref{eq:courant}) on $\mathcal{X}%
(P)\oplus\Omega^{1}(P)$, the symmetric bilinear form is
\[
\langle(\xi_{1},\theta_{1}),(\xi_{2},\theta_{2})\rangle=\frac{1}{2}(i_{\xi
_{1}}\theta_{2}+i_{\xi_{2}}\theta_{1})),
\]
$\rho:TP\oplus T^*P\rightarrow TP$ is projection on the first factor,
and $\mathcal{D}$ is the operator $f\mapsto(0,df)$, so the Leibniz product is
\[
(\xi_{1},\theta_{1})\cdot(\xi_{2},\theta_{2})=([\xi_{1},\xi_{2}],\mathcal{L}%
_{\xi_{1}}\theta_{2}-i_{\xi_{2}}d\theta_{1}).
\]
The Lie algebra $\mathcal{E}/\mathcal{D}(\mathcal{A})$ is thus $\mathcal{X}%
(P)\oplus\Omega^{1}(P)/dC^{\infty}(P)$, on which the bracket, since we can add
to $i_{\xi_{2}}d\theta_{1}$ the exact form $di_{\xi_{2}}\theta_{1}$, is just
the semidirect product bracket of $\mathcal{X}(P)$ acting on $\Omega
^{1}(P)/dC^{\infty}(P)$ by Lie derivation, i.e. $\mathcal{E}/\mathcal{D}%
(\mathcal{A}) = \mathcal{X}(P)\ltimes\Omega^{1}(P)/dC^{\infty}(P)$.

The enveloping Lie algebra is thus a ``double semidirect product''
\[
(\mathcal{X}(P)\ltimes\Omega^{1}(P)/dC^{\infty}(P)) \ltimes(\mathcal{X}%
(P)\times\Omega^{1}(P)).
\]
The action of $(\mathcal{X}(P)\ltimes\Omega^{1}(P)/dC^{\infty}(P)$ on
$\mathcal{X}(P)\times\Omega^{1}(P)$ may be described as follows. Elements of
$\mathcal{X}(P)$ act by Lie derivation on both the vector fields and the
1-forms. A 1-form $\phi$ (modulo exact 1-forms) acts by the nilpotent
operation $(\xi,\theta)\mapsto(0, -i_{\xi}d\phi)$.

\begin{remark}
Pavol ~\v{S}evera has pointed out to us a nice interpretation of the action just
described. First of all, we pass from the Lie algebra to the group which is
the semidirect product of the diffeomorphisms and the abelian group of 1-forms
modulo exact forms. The diffeomorphisms act in the obvious way. To understand
the action of the 1-forms, we think of the action not just on the product
$\mathcal{X}(P)\times\Omega^{1}(P)$, but on the space of subspaces of
$\mathcal{X}(P)\times\Omega^{1}(P)$ which are graphs of 2-forms. Then the
action of a 1-form $\phi$ on a 2-form is simply to add $-d\phi$. These two
operations on 2-forms -- transformation by diffeomorphisms and the addition of
exact forms -- are precisely the operations which may be considered as
symmetries of the variational problem defined by integration of the 2-form
over 2-dimensional submanifolds of $P$.
\end{remark}

\begin{remark}
When $P$ is a compact, oriented manifold of dimension $n$, the space
$\Omega^{1}(P)/dC^{\infty}(P)$ is in natural duality, by integration, with the
space of closed $n-1$--forms on $P$. If $P$ carries a volume element, then the
latter space may be identified with the Lie algebra of volume-preserving
vector fields. The Lie algebra $\mathcal{X}(P)\ltimes\Omega^{1}(P)/dC^{\infty
}(P)$ may then be seen as an enlargement of the Lie algebra of the cotangent
bundle of the group of volume preserving diffeomorphisms. It would be
interesting to relate this interpretation to other aspects of the material in
this paper.
\end{remark}

\begin{remark}
Finally, we note that the constructions in this section can be carried
out equally well in the setting of $(R,\cala)$ $C$-algebras.  These
algebraic objects, introduced in \cite{we:omni}, include as special
cases both the spaces of sections of Courant algebroids and the
omni-Lie algebras of \S \ref{sec-omni}.
\end{remark}

\section{Lie--Yamaguti structures}

\label{sec-yamaguti}

Let $(\mathcal{E},\cdot)$ be a Leibniz algebra, and let $(\mathfrak
{g},\mathfrak{h},f)$ be an enveloping Lie algebra. As in \S2, we have the
reductive decomposition
$\mathfrak{g}\cong\mathfrak{h}\oplus\mathcal{E}_{1/2}$.  We found that
(\ref{eq:proj=original})  recovers the original
skew-symmetrized operation $[\![\cdot,\cdot]\!]$ in $\mathcal{E}$. We will now
show that the skew-symmetric bilinear map $\Delta_{1/2}:\mathcal{E}%
\times\mathcal{E}\rightarrow\mathfrak{h}$ induces an additional operation and
additional structure on $\mathcal{E}$.

More generally, let $\mathfrak{g}$ be a Lie algebra with a reductive
decomposition $\mathfrak{g}=\mathfrak{h}\oplus\mathfrak{m}$, i.e.
$[\mathfrak{h},\mathfrak{h}]\subseteq\mathfrak{h}$ and $[\mathfrak
{h},\mathfrak{m}]\subseteq\mathfrak{m}$. On $\mathfrak{m}$, define bilinear
maps $[\![\cdot,\cdot]\!]:\mathfrak{m}\times\mathfrak{m}\rightarrow
\mathfrak{m}$ and $\Delta:\mathfrak{m}\times\mathfrak{m}\rightarrow
\mathfrak{h}$ by the projections of the Lie bracket:
\begin{align}
\lbrack\![x,y]\!]  &  =\pi_{\mathfrak{m}}\left(  \left[  x,y\right]  \right)
\label{eq:LY-binary}\\
\Delta(x,y)  &  =\pi_{\mathfrak{h}}\left(  \left[  x,y\right]  \right)
\label{eq:LY-delta}%
\end{align}
for $x,y\in\mathfrak{m}$. Then define a ternary product on $\mathfrak{m}$ by
\begin{equation*}
\left\{  x,y,z\right\}  :=[\Delta(x,y),z] 
\end{equation*}
for $x,y,z\in\mathfrak{m}$. It is straightforward to show that $(\mathfrak
{m},[\![\cdot,\cdot]\!],\left\{  \cdot,\cdot,\cdot\right\}  )$ satisfies the
following definition.

\begin{definition}
\label{defn:lie-yamaguti}A \emph{Lie-Yamaguti algebra} $(\mathfrak
{m},[\![\cdot,\cdot]\!],\left\{  \cdot,\cdot,\cdot\right\}  )$ is a vector
space $\mathfrak{m}$ together with a bilinear operation $[\![\cdot
,\cdot]\!]:\mathfrak
{m}\times\mathfrak{m}\to\mathfrak{m}$ and a trilinear operation $\left\{
\cdot,\cdot,\cdot\right\}  :\mathfrak{m}\times\mathfrak
{m}\times\mathfrak{m}\to\mathfrak{m}$ such that, for all $x,y,z,u,v,w\in
\mathfrak{m}$:

\begin{enumerate}
\item [(LY1)]$[\![x,y]\!]=-[\![y,x]\!]$;

\item[(LY2)] $\left\{  x,y,z\right\}  =-\{y,x,z\}$;

\item[(LY3)] $\sum_{(x,y,z)}\left(  [\![[\![x,y]\!],z]\!]+\left\{
x,y,z\right\}  \right)  =0$;

\item[(LY4)] $\sum_{(x,y,z)}\left\{  [\![x,y]\!],z,u\right\}  =0$;

\item[(LY5)] $\left\{  x,y,[\![u,v]\!]\right\}  =[\![\left\{  x,y,u\right\}
,v]\!]+[\![u,\left\{  x,y,v\right\}  ]\!]$;

\item[(LY6)] $\left\{  x,y,\left\{  u,v,w\right\}  \right\}  =\left\{
\left\{  x,y,u\right\}  ,v,w\right\}  +\left\{  u,\left\{  x,y,v\right\}
,w\right\}  +\left\{  u,v,\left\{  x,y,w\right\}  \right\}  $.
\end{enumerate}
\end{definition}

The properties of the binary and ternary operations of a Lie-Yamaguti algebra
can be found in the work of Nomizu \cite{no:invariant} as properties satisfied
by the torsion and curvature tensors, respectively, in a reductive homogeneous
space; we will discuss this further in \S\ref{sec-reductive}. The notion of a
Lie-Yamaguti algebra is a natural abstraction made by K. Yamaguti
\cite{ya:triple}, who called these algebras ``general Lie triple systems''. M.
Kikkawa \cite{ki:geometry} dubbed them ``Lie triple algebras''.

Notice that, if the trilinear product in a Lie-Yamaguti algebra is trivial,
i.e., $\left\{  \cdot,\cdot,\cdot\right\}  \equiv0$, then (LY2), (LY4), (LY5),
and (LY6) are trivial, and (LY1) and (LY3) define a Lie algebra. On the other
hand, if the binary product is trivial, i.e., $[\![\cdot,\cdot]\!]\equiv0$,
then (LY1), (LY4), and (LY5) are trivial, and (LY2), (LY3), and (LY6) define a
Lie triple system.

Now we apply these considerations to the case of a Leibniz algebra
$(\mathcal{E},\cdot)$ with enveloping Lie algebra $(\mathfrak{g},\mathfrak
{h},f)$. The map $\Delta\equiv\Delta_{1/2}:\mathcal{E}\times\mathcal{E}%
\rightarrow\mathfrak{h}$ is given by (\ref{eq:proj-delta}); we repeat it here
for convenience:
\[
\Delta(x,y)=-\frac{1}{4}f\left(  [\![x,y]\!]\right)
\]
for $x,y\in\mathcal{E}$. Therefore the ternary product is given in terms of
the Leibniz and skew-symmetrized Leibniz products by
\begin{equation}
\left\{  x,y,z\right\}  =-\frac{1}{4}[\![x,y]\!]\cdot z
\label{eq:leibniz-ternary}%
\end{equation}
for $x,y,z\in\mathcal{E}$. 
Since $\lambda(x\cdot y + y\cdot x)=0$, we can also write the ternary
product purely in terms of the Leibniz product:
\begin{equation}
\left\{  x,y,z\right\}  =-\frac{1}{4}(x\cdot y)\cdot z.
\label{eq:leibniz-ternary2}%
\end{equation}

Summarizing, we have the following.

\begin{proposition}
\label{prop:leibniz-LY}Let $(\mathcal{E},\cdot)$ be a Leibniz algebra. Then
$(\mathcal{E},[\![\cdot,\cdot]\!],\left\{  \cdot,\cdot,\cdot\right\}  )$ is a
Lie-Yamaguti algebra, where $[\![\cdot,\cdot]\!]$ is the skew-symmetrization
of $\cdot$, and $\left\{  \cdot,\cdot,\cdot\right\}  $ is defined by
(\ref{eq:leibniz-ternary}) or (\ref{eq:leibniz-ternary2}) .
\end{proposition}

\begin{example}
\label{ex:hxV-triple}For the hemisemidirect product $\mathcal{E}=\mathfrak
{h}\ltimes_{H}V$ of Example \ref{ex:hemi  &  demi}, 
a short calculation using (\ref{eq:leibniz-ternary}) and
(\ref{eq:hxV-product2}) shows that the trilinear product in the associated
Lie-Yamaguti algebra $(\mathcal{E},[\![\cdot,\cdot]\!],\left\{  \cdot
,\cdot,\cdot\right\}  )$ is given by
\[
\left\{  (\xi,x),(\eta,y),(\zeta,z)\right\}  =-\frac{1}{4}\left(  [[\xi
,\eta],\zeta],[\xi,\eta]z\right)
\]
for $\xi,\eta,\zeta\in\mathfrak{h}$, $x,y,z\in V$.
\end{example}

\begin{example}
\label{ex:courant-triple} 
For the Courant bracket (\ref{eq:courant}), the associated trilinear
product is
\begin{equation*}
\{(\xi_{1},\theta_{1}),(\xi_{2},\theta_{2}),(\xi_{3},\theta_{3})\}
=-\frac{1}{4}([[\xi_{1},\xi_{2}],\xi_3],\mathcal{L}_{[\xi_1,\xi_2]}\theta_{3}-i_{\xi_3}d 
(\mathcal{L}_{\xi_1}\theta_2 - \mathcal{L}_{\xi_2}\theta_1)).
\end{equation*}

Another approach to this bracket, and to spaces of sections of more
general Courant algebroids, is to consider them as \emph{strongly homotopy
Lie algebras}   \cite{ro:courant,ro-we:courant}. Here, the
not-quite-Lie algebra also carries a differential $d$, and the
Jacobiator $\sum_{(x,y,z)} [\![[\![x,y]\!],z]\!]$ is expressed
as the differential applied to a completely antisymmetric trilinear
product.  A higher-order Jacobiator of the trilinear product is again
a differential, and so on.  In the Lie-Yamaguti approach, it is
essential that the trilinear operation \emph{not} be completely
antisymmetric.  
\end{example}

We showed that a Lie algebra with a reductive decomposition naturally induces
the structure of a Lie-Yamaguti algebra on the reductive complement to the
subalgebra. Now we show the converse: for any Lie-Yamaguti algebra
$\mathfrak{m}$, there exists a Lie algebra $\mathfrak{g}$ with a reductive
decomposition $\mathfrak{g}=\mathfrak{h}\oplus\mathfrak{m}$ such that the
induced Lie-Yamaguti structure on $\mathfrak{m}$ agrees with the original one.

Let $(\mathfrak{m},[\![\cdot,\cdot]\!],\left\{  \cdot,\cdot,\cdot\right\}  )$
be a Lie-Yamaguti algebra, and let $\mathrm{Der}(\mathfrak{m})$ denote the Lie
subalgebra of $\mathfrak
{gl}(\mathfrak{m})$ consisting of derivations of both the bilinear and
trilinear products. Let $\mathfrak{h}$ be a Lie algebra with a derivation
action $\mathfrak{h}\rightarrow\mathrm{Der}(\mathfrak{m})$, and let
$\Delta:\mathfrak{m}\times\mathfrak
{m}\rightarrow\mathfrak{h}$ be a skew-symmetric bilinear mapping satisfying
the following properties:%
\begin{gather}
\Delta(x,y)z=\left\{  x,y,z\right\}  \label{eq:delta1}\\
\left[  \xi,\Delta(x,y)\right]  =\Delta(\xi x,y)+\Delta(x,\xi
y)\label{eq:delta2}\\
\Delta\left(  [\![x,y]\!],z\right)  +\Delta\left(  [\![y,z]\!],x\right)
+\Delta\left(  [\![z,x]\!],y\right)  =0\label{eq:delta3}%
\end{gather}
for all $x,y,z\in\mathfrak{m}$, $\xi\in\mathfrak{h}$. On the vector space
$\mathfrak{g}=\mathfrak{h}\oplus\mathfrak{m}$, define a skew-symmetric bracket
$[\cdot,\cdot]:\mathfrak{g}\times\mathfrak{g}\rightarrow\mathfrak{g}$ by
\begin{equation}
\left[  \xi+x,\eta+y\right]  =\left(  \left[  \xi,\eta\right]  +\Delta
(x,y)\right)  +\left(  \xi y-\eta x+[\![x,y]\!]\right)  .\label{eq:g-bracket}%
\end{equation}

\begin{proposition}
$\mathfrak{g}$ is a Lie algebra, and $\mathfrak{g}=\mathfrak{h}\oplus
\mathfrak{m}$ is a reductive decomposition. The Lie-Yamaguti bilinear and
trilinear products on $\mathfrak{m}$ induced by the decomposition agree with
$[\![\cdot,\cdot]\!]$ and $\left\{  \cdot,\cdot,\cdot\right\}  $, respectively.
\end{proposition}

\begin{proof}
It is straightforward to check that (\ref{eq:delta1})-(\ref{eq:delta3}) and
(LY3) give the Jacobi identity for the bracket (\ref{eq:g-bracket}). That
$\left[  \mathfrak{h},\mathfrak{m}\right]  \subseteq\mathfrak{m}$ is immediate
from (\ref{eq:g-bracket}), and so the decomposition is reductive. From
(\ref{eq:g-bracket}) we have $\pi_{\mathfrak{m}}\left(  \left[  x,y\right]
\right)  =[\![x,y]\!]$ and $\left[  \pi_{\mathfrak{h}}\left(  \left[
x,y\right]  \right)  ,z\right]  =\left[  \Delta\left(  x,y\right)  ,z\right]
=\left\{  x,y,z\right\}  $ for $x,y,z\in\mathfrak{m}$, which proves the
remaining assertion.
\end{proof}

\begin{definition}
Let $(\mathfrak{m},[\![\cdot,\cdot]\!],\left\{  \cdot,\cdot,\cdot\right\}  )$
be a Lie-Yamaguti algebra, let $\mathfrak{h}$ be a Lie algebra with a
derivation action $\mathfrak{h}\rightarrow\mathrm{Der}(\mathfrak{m})$, and let
$\Delta:\mathfrak{m}\times\mathfrak{m}\rightarrow\mathfrak{h}$ be a
skew-symmetric bilinear mapping satisfying (\ref{eq:delta1})-(\ref{eq:delta3}%
). The Lie algebra $(\mathfrak{g},[\cdot,\cdot])$, where $\mathfrak
{g}=\mathfrak{h}\oplus\mathfrak{m}$ and $[\cdot,\cdot]$ is given by
(\ref{eq:g-bracket}), is called an \emph{enveloping Lie algebra} of
$(\mathfrak{m},[\![\cdot,\cdot]\!],\left\{  \cdot,\cdot,\cdot\right\}  )$.
\end{definition}

Since a Leibniz algebra has a natural Lie-Yamaguti structure, the coincidence
of their notions of enveloping Lie algebra is both clear and expected.

\begin{proposition}
Let $(\mathcal{E},\cdot)$ be a Leibniz algebra with enveloping Lie algebra
$(\mathfrak{g},\mathfrak{h},f)$. Then $\mathfrak{g}\cong\mathfrak{h}%
\oplus\mathcal{E}_{1/2}$ is an enveloping Lie algebra of the induced
Lie-Yamaguti algebra $(\mathcal{E},[\![\cdot,\cdot]\!],\left\{  \cdot
,\cdot,\cdot\right\}  )$.
\end{proposition}

To conclude this section, we follow Yamaguti \cite{ya:triple} to show that
\emph{every} Lie-Yamaguti algebra has an enveloping Lie algebra.

Let $(\mathfrak{m},[\![\cdot,\cdot]\!],\left\{  \cdot,\cdot,\cdot\right\}  )$
be a Lie-Yamaguti algebra. Define a bilinear mapping $\delta:\mathfrak
{m}\times\mathfrak{m}\to\mathfrak
{gl}(\mathfrak{m})$ by
\begin{equation}
\delta(x,y)z=\left\{  x,y,z\right\}  \label{eq:delta}%
\end{equation}
for $x,y,z\in\mathfrak{m}$. By (LY2), $\delta$ is skew-symmetric. By (LY5) and
(LY6), $\delta(x,y)\in\mathrm{Der}(\mathfrak{m})$ for all $x,y\in\mathfrak{m}%
$. We call $\delta(x,y)$ an \emph{inner derivation} of $\mathfrak{m}$. Let
$\mathrm{IDer}(\mathfrak{m})$ denote the subspace of $\mathrm{Der}%
(\mathfrak{m})$ spanned by the inner derivations. By (LY6), $\mathrm{IDer}%
(\mathfrak{m})$ is a\ Lie algebra. Now let $\mathfrak{h}$ be any Lie algebra
satisfying $\mathrm{IDer}(\mathfrak{m})\subseteq\mathfrak{h}\subseteq
\mathrm{Der}(\mathfrak{m})$. Then (\ref{eq:delta}) implies (\ref{eq:delta1})
and (\ref{eq:delta2}) directly, while (\ref{eq:delta3}) follows from both
(\ref{eq:delta}) and (LY4).

Summarizing, we have established the following \cite{ya:triple}.

\begin{proposition}
\label{prop:LY-envelope}Every Lie-Yamaguti algebra has an enveloping Lie algebra.
\end{proposition}

\begin{remark}
It was essentially by this route, after guessing the ternary product
(\ref{eq:leibniz-ternary}), that we first discovered the enveloping Lie
algebras of Leibniz algebras.
\end{remark}


\begin{remark}
It would be interesting to find conditions on a reductive decomposition
$\mathfrak{g}=\mathfrak{h}\oplus\mathfrak{m}$ (equivalently, on a Lie-Yamaguti
algebra $\mathfrak{m}$) which would insure that $\mathfrak{m}$ is the
skew-symmetrization of a Leibniz algebra. What is necessary, of course, is
that $\Delta:\mathfrak{m}\times\mathfrak{m}\rightarrow\mathfrak{h}$ factors
into an $\mathfrak{h}$-equivariant linear map $-\frac{1}{4}f:\mathfrak
{m}\rightarrow\mathfrak{h}$ and the binary product $[\![\cdot,\cdot
]\!]:\mathfrak{m}\times\mathfrak{m}\rightarrow\mathfrak{m}$. The question then
becomes: what conditions on the decomposition guarantee the existence of such
a factorization?
\end{remark}

\section{Reductive homogeneous spaces and loops}

\label{sec-reductive}
Let $G$ be a Lie group with Lie algebra $\mathfrak{g}$,
let $H\subseteq G$ be a closed subgroup with Lie algebra $\mathfrak
{h}\subseteq\mathfrak
{g}$, and let $M=G/H$. The homogeneous space $M$ is said to be
\emph{reductive} if there exists a reductive complement $\mathfrak{m}$ of
$\mathfrak{h}$ in $\mathfrak{g}$, i.e.,
\begin{gather}
\mathfrak{g}=\mathfrak{h}\oplus\mathfrak{m}\label{eq:reductive1}\\
\mathrm{Ad}_{G}(H)\mathfrak{m}\subseteq\mathfrak{m} \label{eq:reductive2}%
\end{gather}
Condition (\ref{eq:reductive2}) implies that
\begin{equation}
\mathrm{ad}_{\mathfrak{g}}(\mathfrak{h})\mathfrak{m}\subseteq\mathfrak{m}.
\label{eq:reductive3}%
\end{equation}
When $H$ is connected, (\ref{eq:reductive3}) implies (\ref{eq:reductive2}). We
may identify $\mathfrak{m}$ with the tangent space $T_{\pi(e)}(M)$ where $e\in
G$ is the identity element and $\pi:G\to M$ is the canonical projection.

Now let $(\mathcal{E},\cdot)$ be a Leibniz algebra with enveloping Lie algebra
$(\mathfrak{g},\mathfrak{h},f)$, and let $H$ be a Lie group with Lie algebra
$\mathfrak{h}$. Suppose that the derivation action of $\mathfrak{h}$ on
$\mathcal{E}$ lifts to an automorphism action $H\rightarrow\mathrm{Aut}%
(\mathcal{E})$, i.e. $h(x\cdot y)=hx\cdot hy$ for $h\in H$, $x,y\in
\mathcal{E}$. Let $G=H\ltimes\mathcal{E}$ be the semidirect product group,
$\mathcal{E}$ considered as usual as an abelian Lie group. Then $G$ is a Lie
group with Lie algebra $\mathfrak{g}=\mathfrak{h}\ltimes\mathcal{E}$.

For $s\neq 0$, we already know that (\ref{eq:reductive3}) is satisfied with
$\mathfrak{m}=\sigma_s(\mathcal{E})$. Now since $H$ acts by automorphisms,
the mapping $f:\mathcal{E}\rightarrow\mathfrak{h}$ is $H$-equivariant:
for all $h\in H$, $x\in\mathcal{E}$,
\begin{equation}
\mathrm{Ad}(h)f(x)=f(hx).\label{eq:f-H-equivariant}
\end{equation}
It follows that the sections
$\sigma_s :\mathcal{E}\to\mathfrak{g}$ defined by
$ x\mapsto (sf(x),x)$ are also
$H$-equivariant, i.e.
\[
h\sigma_s(x) = h(sf(x),x) =(sf(hx),hx)=\sigma_s(hx).
\]
Thus (\ref{eq:reductive2}) holds. Therefore $G/H$ is a reductive
homogeneous space.

The vector space projection
$\pi_{\mathcal{E}}:\mathfrak{g}\to\mathcal{E}$ defined by $(\xi,x)\mapsto x$
exponentiates to the group projection
 $\hat{\pi}_{\mathcal{E}}: G\to\mathcal{E}; (h,x)\mapsto x$.
The sections
$\sigma_s$
exponentiate to the group sections
$\hat{\sigma}_s :\mathcal{E}\to G$ defined by
\begin{equation*}
\hat{\sigma}_s(x) = (\exp (sf(x)),x) 
\end{equation*}
For $(h,x)\in G=H\ltimes\mathcal{E}$, we have
\[
(h,x) = \hat{\sigma}_s(x)( \exp (-sf(x)) h,0).
\]
This is clearly a unique factorization of $(h,x)$ into an element of
the image
$\hat{\cale}_s  =  \hat{\sigma}_s(\mathcal{E})$ and an element of
$H\cong H\times \{ 0\}$. This implies that $\hat{\cale}_s $
is a \emph{left transversal} of $H$ in $G$, i.e. a subset of $G$
consisting of one representative of each left coset in $G/H$.

Summarizing, we use the semidirect product structure of $G$ to
identify the homogeneous space $(H \ltimes \mathcal{E})/H$ with $\mathcal{E}$ itself,
and we have a distinguished family of sections
$\hat{\sigma}_s : \mathcal{E}\to G$ whose images are transversals
of the subgroup $H$.

In order to motivate our next construction for $\mathcal{E}$,
we will temporarily forget about differentiable structure. Let $G$
be a group acting transitively on a set $X$. Fix a distinguished
element $e\in X$, and let $H$ be the isotropy subgroup of $e$.
Associated to $e$ is the canonical projection
$\pi_X: G\to X;g\mapsto g(e)$. Assume
that $\phi:X\to G$ is a section
of $\pi_X$, i.e. $\pi_X ( \phi(x)) = \phi(x)(e) = x$ for all
$x\in X$, and assume that $\phi(e) = 1\in G$.
Let $\pi_H : G\to H$ be the corresponding projection
onto $H$ defined by $\pi_H(g) = (\phi (\pi_X(g)))^{-1} g$.
We define a binary operation $\diamond : X\times X\to X$ and
a mapping $l : X\times X\to H$ by
\begin{align}
x\diamond y &= \pi_X (\phi(x)\phi(y)) \label{eq:loop-op}\\
l(x,y) &= \pi_H (\phi(x)\phi(y)) \label{eq:trans-map}
\end{align}
for $x,y\in X$. Then $(X,\diamond)$ is a \emph{left loop} \cite{ki:loops},
i.e., given $a,b\in X$, the equation $a\diamond x=b$ has a unique solution
$x\in X$, and $e\diamond a=a\diamond e=a$ for all $a\in X$. (Conversely,
every left loop can be realized in this way, specifically, as a left
transversal in a group.) For $a,b\in X$, the action of $l(a,b)\in H$
on $X$ defines a permutation $L(a,b):X\to X$, called a
\emph{left inner mapping}, by $L(a,b)(c) = l(a,b)(c)$ for all $c\in X$.
Left inner mappings measure the nonassociativity of $(X,\diamond)$;
they are equivalently defined by the equation
$a\diamond(b\diamond c)=(a\diamond b)\diamond L(a,b)(c)$ for
$a,b,c\in X$.

Two useful conditions on sections $\phi:X\to G$ are the following:

\begin{enumerate}
\item[(H1)] For each $x\in X$, there exists (a necessarily unique)
$x'\in X$ such that
\begin{equation*}
\phi(x)^{-1} = \phi(x') .
\end{equation*}

\item[(H2)] For all $x\in X$, $h\in H$,
\begin{equation*}
h\phi(x) h^{-1} = \phi(h(x)).
\end{equation*}
\end{enumerate}
If (H1) holds, then the left loop $(X,\diamond)$ satisfies
the \emph{left inverse property}, which means that
$a'\diamond(a\diamond b)=a\diamond(a'\diamond b)=b$ for all
$a,b\in X$. If (H2) holds, then action of $H$ on $X$ is
by automorphisms of $(X,\diamond)$, i.e.
\[
h(x\diamond y) = h(x)\diamond h(y)
\]
for all $x,y\in X$, $h\in H$. In particular, every left inner mapping
$L(a,b)$ is an automorphism, and in this case, $(X,\diamond)$ is said
to have the $A_l$ (or \emph{left} $A$ or \emph{left special})
\emph{property}. A left loop with both the left inverse and $A_{l}$
properties is said to be \emph{homogeneous}. For more on these matters,
see \cite{ki:loops} and the references therein.

\begin{remark}
Let $G$ be a Lie group with closed (Lie) subgroup $H\subseteq G$ and
a smooth section $\phi:G/H \to G$ of the natural projection
$G\to G/H$. Kikkawa \cite{ki:geometry} showed that if (H1) and (H2)
hold, then $G/H$ is a reductive homogeneous space. The converse
problem, which is to characterize reductive homogeneous spaces
$G/H$ such that there exists a smooth section $\phi:G/H \to G$
satisfying (H1) and (H2), is still open.
\end{remark}

We apply the preceding notions to the sections
$\hat{\sigma}_s:\mathcal{E}\to G$ in the group $G=H\ltimes\mathcal{E}$.
For $x,y\in\mathcal{E}$, we have
\begin{align*}
\hat{\sigma}_s(x) \hat{\sigma}_s(y) &=
(\exp(sf(x)),x)(\exp(sf(y)),y)\\
&  =(\exp(sf(x))\exp(sf(y)),x+\exp(sf(x))y).
\end{align*}
Following (\ref{eq:loop-op}), we take the projection of this
product onto $\mathcal{E}$ to define a family of left loop
structures $(\mathcal{E},\diamond_s)$ by
\begin{equation}
x\diamond_s y=x+\exp(sf(x))y \label{eq:E-loop}
\end{equation}
for $x,y\in\mathcal{E}$.  Also, following (\ref{eq:trans-map}), we
define a corresponding family of maps
$l_s:\mathcal{E}\times \mathcal{E}\to H$ by
\begin{equation*} 
l_s(x,y)=  \exp(s(-x-\exp(sf(x))y))\exp(sf(x))\exp(sf(y))
\end{equation*}
for $x,y\in\mathcal{E}$.

Before considering homogeneity, we first note the following
important property.

\begin{proposition}\label{prop:independent}
The left loops $(\mathcal{E},\diamond_s)$ defined by
(\ref{eq:E-loop}) depend only on the Leibniz algebra
structure $(\mathcal{E},\cdot )$, and not on the choice
of enveloping algebra $(\mathfrak{g},\mathfrak{h},f)$.
In particular,
\begin{equation}
x\diamond_s y = x + \exp(s\lambda(x))y \label{eq:E-loop2}
\end{equation}
for all $x,y\in\mathcal{E}$.
\end{proposition}

\begin{proof} Indeed, for $x,y\in\mathcal{E}$, we have
\[
\exp(sf(x))y = \exp(s\lambda(x))y
\]
using (2.13). This establishes the equivalence of
(\ref{eq:E-loop}) and (\ref{eq:E-loop2}).
\end{proof}

Next we turn to homogeneity. For each $x\in\mathcal{E}$, we have
\begin{align}
\hat{\sigma}_s(x)^{-1} &= (\exp(sf(x)),x)^{-1} \nonumber \\
&= (\exp(-sf(x)), -\exp(-sf(x))x). \label{eq:inv-tmp}
\end{align}
Now
\begin{align*}
\exp(-sf(x)) &= \exp (-s \mathrm{Ad}(\exp(-sf(x))) f(x)) \\
&= \exp (s f (-\exp (-sf(x)) x)),
\end{align*}
using (\ref{eq:f-H-equivariant}). Applying this to
(\ref{eq:inv-tmp}), we have
\[
\hat{\sigma}_s(x)^{-1} = \hat{\sigma}_s(-\exp(-sf(x))x) .
\]
Thus (H1) holds, and the left
loops $(\mathcal{E},\diamond_s)$ satisfy the left inverse
property. In particular, the inverse of $x$ in
$(\mathcal{E},\diamond_s)$ is given by
\[
x^{\prime} = -\exp(-s \lambda(x)) x .
\]

The $\mathrm{Ad}_{G}(H)$-invariance of $\mathcal{E}$
implies that $H$ normalizes $\hat{\sigma}_s(\mathcal{E})$.
More precisely, for $h\in H$, $x\in\mathcal{E}$,
\begin{align*}
(h,0)\hat{\sigma}_s(x)(h^{-1},0)
&= (h\exp(sf(x))h^{-1},hx)\\
&= (\exp(s\mathrm{Ad}(h)f(x)),hx)\\
&= (\exp(sf(hx),hx)\\
&= \hat{\sigma}_s(hx).
\end{align*}
Thus (H2) holds, and so $(\mathcal{E},\diamond_s)$ has the
$A_l$ property. Summarizing, we have the following.

\begin{proposition}\label{prop:left-loop} Let $(\mathcal{E},\cdot)$
be a Leibniz algebra, and for $s\neq 0$, let
$\diamond_s:\mathcal{E}\times \mathcal{E}\to \mathcal{E}$ be
defined by
\[
x\diamond_s y = x + \exp(s\lambda(x))y
\]
for $x,y\in\mathcal{E}$. Then $(\mathcal{E},\diamond_s)$ is a
homogeneous left loop for which the binary bracket in the associated
Lie-Yamaguti algebra is $2s$ times the skew-symmetrized Leibniz
product.  In particular, the skew-symmetrized Leibniz
product itself is recovered when $s=1/2$.
\end{proposition}

In addition, $(\mathcal{E},\diamond_s)$ 
is a \emph{geodesic} left loop, which
means that it agrees with the natural local left loop structure defined in a
neighborhood of $0\in\mathcal{E}$ by parallel transport of geodesics along
geodesics \cite{ki:on,ki:geometry}.  We note that the loop structure
$\diamond_1$ was described in the case of Lie algebras by Kikkawa in
Proposition 4 of \cite{ki:projectivity}.

\begin{example}
\label{ex:omniloop}
For the hemisemidirect product Leibniz algebra  $\frakh \ltimes_H V$ of
Example \ref{ex:hemi  & demi}, built from the
representation $\rho : \frakh \rightarrow \mathfrak{gl}(V)$,
 the left loop structures coming from 
Proposition \ref{prop:left-loop} are 
\begin{equation*}
(\xi,x)\diamond_s (\eta,y) = (\xi +\exp(\mathrm{ad}(s\xi))(\eta),x+\exp
(\rho (s\xi))y).
\end{equation*}
\end{example}

\begin{example}
\label{ex:courantloop}
For the Courant bracket (\ref{eq:courant}) on $\mathcal{X}%
(P)\oplus\Omega^{1}(P)$,  the left loop structures coming from
Proposition \ref{prop:left-loop} are
\[
(\xi_1,\theta_1)\diamond_s (\xi_2,\theta_2)=(\xi,\theta),\]
where
\[
\xi=\xi_1+(\exp s\xi_1)^*\xi_2,\]
and
\[\theta = \theta_1+(\exp s\xi_1)^* \theta_2- (\exp s\xi_1)^*
(\xi_2)\backl d \int_0^s (\exp t\xi_1)^* \theta_1 dt. \]

\end{example}

\section{Connections}

\label{sec-connections}

 The
following classification theorem for $G$-invariant connections on reductive
homogeneous spaces is due to Nomizu (\cite{no:invariant}, Thm. 8.1).

\begin{proposition}
\label{prop:nomizu}Let $M=G/H$ be a reductive homogeneous space with fixed
decomposition (\ref{eq:reductive1})-(\ref{eq:reductive2}). There exists a
one-to-one correspondence between the set of all $G$-invariant connections on
$M$ and the set of all $\mathrm{Ad}_{G}(H)$-equivariant bilinear mappings
$\alpha:\mathfrak{m}\times\mathfrak{m}\to\mathfrak{m}$.
\end{proposition}

The $\mathrm{Ad}_{G}(H)$-equivariance of $\alpha:\mathfrak{m}\times
\mathfrak{m}\to\mathfrak{m}$ means that $\mathrm{Ad}(h)\alpha(X,Y)=\alpha
(\mathrm{Ad}(h)X,\mathrm{Ad}(h)Y)$ for all $h\in H$, $X,Y\in\mathfrak{m}$.
Equivalently, if one thinks of $(\mathfrak{m},\alpha)$ as a nonassociative
algebra, $\mathrm{Ad}_{G}(H)$ is a subgroup of the automorphism group
$\mathrm{Aut}(\mathfrak{m},\alpha)$.

Let $(\mathfrak{m},[\![\cdot,\cdot]\!],\left\{  \cdot,\cdot,\cdot\right\}  )$
be the Lie-Yamaguti algebra structure on $\mathfrak{m}$ determined by the
decomposition (\ref{eq:reductive1})-(\ref{eq:reductive2}); thus $[\![\cdot
,\cdot]\!]=[\cdot,\cdot]_{\mathfrak{m}}$ and $\left\{  \cdot,\cdot
,\cdot\right\}  =[[\cdot,\cdot]_{\mathfrak{h}},\!]$, where, as before, the
subscripts indicate projections. 
The \emph{zero} bilinear mapping on $\mathfrak{m}$ corresponds to the 
\emph{canonical connection} (of the $2$nd kind) on $M$. This connection 
is characterized by
the following property: for each $X\in\mathfrak{m}$, parallel displacement of
tangent vectors along the curve $\pi(\exp tX)$ is the same as the translation
of tangent vectors by the natural action of $\exp tX$ on $M$. The torsion and
curvature are
\begin{align}
T(X,Y)  &  =-[\![X,Y]\!]\label{eq:canon-tor}\\
R(X,Y)Z  &  =-\left\{  X,Y,Z\right\}  . \label{eq:canon-curv}%
\end{align}
for $X,Y,Z\in\mathfrak{m}$.


If $M=G/H$ is a homogeneous left loop, Kikkawa showed that the canonical
connection can be constructed directly from the loop multiplication
(\cite{ki:geometry}, Thm. 3.7). The torsion and curvature tensors then
define a Lie-Yamaguti algebra structure on $\mathfrak{m}$ by
(\ref{eq:canon-tor})-(\ref{eq:canon-curv}), which is considered to
be the tangent algebra structure of the homogeneous loop. In case
$H=\{1\}$, the canonical connection is just Cartan's $(-)$-connection
on $G$, and the corresponding Lie-Yamaguti algebra is the Lie algebra
of $G$ (\cite{ki:geometry}, Ex. 3.3).

Now let $(\mathcal{E},\cdot)$ be a Leibniz algebra with enveloping Lie algebra
$\mathfrak{g}=\mathfrak{h}\oplus \mathcal{E}_s$, and homogeneous left
loop structure $(\mathcal{E},\diamond_s)$. We identify $\mathcal{E}$
with the reductive homogeneous space $G/H$. Following Kikkawa
\cite{ki:geometry} (see also Miheev and Sabinin \cite{mi:quasigroups}),
one finds that the canonical connection is given by 
\begin{equation}\label{eq:connection}
(\nabla_X Y)(x) = DY(x)X(x) - s X(x)\cdot Y(x)
\end{equation}
for $x\in \mathcal{E}$. Here $X$ and $Y$ are vector fields on
$\mathcal{E}$ which we are identifying with mappings $X,Y:\mathcal{E}
\to \mathcal{E}$. 

We observe that: (1) the reductive decompositions
$\mathfrak{g} = \mathfrak{h} \oplus \mathcal{E}_s$ give a
one-parameter deformation of the semidirect product
$\mathfrak{g} = \mathfrak{h}\ltimes \mathcal{E}$; (2) the connection
 (\ref{eq:connection}) is
 a  deformation of the standard flat,
torsion-free connection on $\mathcal{E}$; (3) the loop structures
$\diamond_s$ form a deformation of the addition operation on $\cale$.

\section{Further Questions}
\label{sec-questions}

We recall from \S\ref{sec-courant} that
 Courant's algebra $\mathcal{E} = \mathcal{X}(P)\oplus
\Omega^{1}(P)$ is a Courant algebroid; i.e. it is also 
the $C^{\infty}$-module of sections of a vector bundle
over $P$, and the bracket satisfies identities which relate it to the module
structure. A simpler version of these identities defines \emph{Lie
algebroids} \cite{ma:lie}, 
which are the infinitesimal objects associated to Lie groupoids.
The sections of a Lie algebroid form a Lie algebra which acts by derivations
of $C^{\infty}(P)$, and the corresponding group of ``bisections'' of the
groupoid acts by automorphisms, i.e. by diffeomorphisms of $P$. It has been
our hope to find a group-like object associated to $\mathcal{E} =
\mathcal{X}(P)\oplus\Omega^{1}(P)$ which has something like an action on $P$,
and which can be considered as the sections of some kind of nonassociative
generalization of a groupoid (a loopoid?). So far, we have not succeeded.
The difficulty might be related to the absence of a 
natural adjoint representation of a Lie algebroid on itself (as opposed to the
adjoint representation of the Lie algebra of sections). Perhaps a
weak version of the adjoint representation, such as is described in the
appendices of \cite{ev-lu-we:transverse}, could be a model for what we seek in
the case of Courant algebroids.

Finally, we are left with the problem of constructing a group-like
object attached to a Leibniz algebra in such a way that the object is
a group when the Leibniz algebra is a Lie algebra.  A possible
approach to this problem is via path spaces.  At the end of Chapter 1 in
\cite{du-ko:lie}, Duistermaat and Kolk prove ``Lie's third theorem''
by beginning with a Lie algebra $\cale$ and defining a Banach Lie group structure
on the space $\calp(\cale)$ of continuous paths $\gamma:[0,1]\to\cale$.  When
the Lie algebra of this group is identified with $\calp(\cale)$
itself, the integration map $I:\gamma\mapsto \int_0^1 \gamma(t)dt$
 is found to be a homomorphism
from $\calp(\cale)$ to $\cale$
.  The closed ideal $\ker I$ integrates to a
normal Lie subgroup $\calp_0(\cale) \subset \calp(\cale)$ which is
shown\footnote{Nothing comes for free.  The proof that
$\calp_0(\cale)$ is closed relies on the same vanishing theorem for
the second cohomology of a finite dimensional simply connected Lie
group which goes into other proofs of Lie III.  Of course, the result
holds only when $\cale$ is finite dimensional, as it should.} to be
closed.  The quotient $\calp(\cale)/ \calp_0(\cale)$ is then a Lie
group whose Lie algebra is isomorphic to $\cale$.  

Cattaneo and Felder \cite{ca-fe:poisson} have used a similar path
space construction to construct symplectic groupoids from Poisson
manifolds.  Their construction is a variation of the Duistermaat--Kolk
construction applied to the cotangent bundle Lie algebroid.  Rather
than using an associative product corresponding to pointwise multiplication
of group(oid) paths, they use concatenation of Lie algebroid paths, which
becomes a groupoid structure only after an equivalence relation is
applied; the idea should be extendible to arbitrary Lie algebroids.  Here,
the resulting groupoid may have singularities; only the local groupoid
is smooth.

We have begun to investigate the Duistermaat-Kolk construction when
$(\cale,\cdot)$ is a Leibniz algebra.  If we use the same formula as
in \cite{du-ko:lie}, the multiplication on $\calp(\cale)$ is,
remarkably, still associative, so we still have a Banach Lie group.
On the other hand, the kernel $\calp_0(\cale)$ of the integration map
is no longer an ideal; in fact it is not even a subalgebra unless
$\cale$ is a Lie algebra.  We can still recover the skew-symmetrized
Leibniz bracket by identifying $\cale$ with the constant paths, and
projecting the path space Lie bracket along $\calp_0(\cale)$, but
since the latter is not a Lie subalgebra, it is not clear how to pass
the group product to a quotient space.

Although the kernel $\calp_0(\cale)$ is not a subalgebra, there is a
different complement of the constant paths which \emph{is} a
subalgebra, namely the kernel of the evaluation map
$E_{1/2}:\gamma\mapsto \gamma(1/2)$ (a crude ``midpoint
approximation'' to $I$).  The kernel $\ker E_{1/2}$ is also the
corresponding Lie subgroup of $\calp(\cale)$.  Although the adjoint
action of $\ker E_{1/2}$ does not leave the constant paths invariant,
we can still use projection along $\ker E_{1/2}$ to construct Lie
algebra and loop structures on $\cale$.  The projected bracket turns
out to be the antisymmetrized Leibniz product $\bleft~.~\bright$,
while the projected product is the same multiplication
$\diamond_{1/2}$ which was obtained in \S\ref{sec-reductive} by a
completely different method.

In a sense, then, we are back where we started.  What we still need is
a construction which incorporates the best properties of $I$ (which
produces the right group when $\cale$ is a Lie algebra) and $E_{1/2}$
(whose kernel is a subalgebra for any Leibniz algebra $\cale$), to reduce
$\calp(\cale)$ to manageable size.  Our search will continue.


\begin{thebibliography}{99}

\bibitem{ca-fe:poisson}Cattaneo, A.S., and Felder, G., Poisson sigma models and
symplectic groupoids, preprint math.SG/0003023.

\bibitem{co:dirac}Courant, T.J., Dirac manifolds, \emph{Trans. A.M.S.}
\textbf{319} (1990), 631--661.

\bibitem{du-ko:lie}
Duistermaat, J.J., and Kolk, J.A.C., {\em Lie Groups},
Berlin-Heidelberg-New York, Springer-Verlag, 1999.

\bibitem{ev-lu-we:transverse}Evens, S., Lu, J.-H., and Weinstein, A.,
Transverse measures, the modular class, and a cohomology pairing for Lie
algebroids, \emph{ Quart. J. Math} \textbf{50} (1999), 417--436.

\bibitem{ki:on}Kikkawa, M., On local loops in affine manifolds, \emph{J. Sci.
Hiroshima Univ. Ser. A-I} \textbf{28} (1964), 199-207.

\bibitem{ki:geometry}Kikkawa, M., Geometry of homogeneous Lie loops,
\emph{Hiroshima Math. J.} \textbf{5} (1975), 141--179.

\bibitem{ki:projectivity}Kikkawa, M., Projectivity of left loops on
$R^n$, \emph{Mem. Fac. Sci. Shimane Univ.} \textbf{22} (1988), 33-41.

\bibitem{ki:loops}Kinyon, M.K., and Jones, O., Loops and semidirect products,
\emph{Comm. Algebra}, to appear.

\bibitem{ko:from}Kosmann-Schwarzbach, Y., From Poisson algebras to
Gerstenhaber algebras, \emph{Ann. Inst. Fourier (Grenoble)} \textbf{46}
(1996), 1243--1274.

\bibitem{li-we-xu:manin}Liu, Z.-J., Weinstein, A., and Xu, P., Manin triples
for Lie bialgebroids, \emph{J. Diff. Geom} \textbf{45} (1997), 547--574.

\bibitem{lo:version}Loday, J.L., Une version non commutative des alg\`{e}bres
de Lie: les alg\`{e}bres de Leibniz, \emph{Enseign. Math.} \textbf{39}
(1993), 269--293.

\bibitem{lo:tensor}Loday, J.L., and Pirashvili, T., The tensor category of
linear maps and Leibniz algebras, \emph{Georgian Math. J.} \textbf{5}
(1998), 263-276.

\bibitem{ma:lie} 
Mackenzie, K., {\em Lie Groupoids and Lie Algebroids in Differential
Geometry}, LMS Lecture Notes Series, {\bf 124}, Cambridge Univ. Press, 1987.



\bibitem{mi:quasigroups} Miheev, P.O.,  and Sabinin, L.V.,  Quasigroups
and differential geometry, Chapter XII in \emph{Quasigroups and Loops:
Theory and Applications}, Chein, O., Pflugfelder, H.O.,  and 
Smith, J.D.H.  (eds.), Berlin: Heldermann-Verlag, 1990.

\bibitem{myung:malcev}Myung, H.C., \emph{Malcev-Admissible Algebras},
Boston-Basel-Stuttgart: Birkh\"{a}user, 1986.

\bibitem{no:invariant}Nomizu, K., Invariant affine connections on homogeneous
spaces, \emph{Amer. J. Math.} \textbf{76} (1954), 33-65.

\bibitem{ro:courant}Roytenberg, D., Courant algebroids, derived brackets and
even symplectic supermanifolds, Ph.D. thesis, University of California,
Berkeley (1999), preprint math.DG/9910078.

\bibitem{ro-we:courant}Roytenberg, D., and Weinstein, A., Courant algebroids
and strongly homotopy Lie algebras, \emph{Lett. Math. Phys.} \textbf{46}
(1998), 81--93.

\bibitem{we:omni}Weinstein, A., Omni-Lie algebras, preprint math.RT/9912190,
to appear in \emph{RIMS K\^{o}ky\^{u}roku}.

\bibitem{ya:triple}Yamaguti, K., On the Lie triple system and its
generalization, \emph{J. Sci. Hiroshima Univ. Ser. A} \textbf{21} (1957/1958), 155--160.
\end{thebibliography}
\end{document}